\documentclass[11pt]{article}
\usepackage{latexsym,amssymb,amsmath,amsfonts,amsthm}
\usepackage{graphics,graphicx,mathrsfs,subfigure}
\usepackage{color,xspace}
\newcommand{\be}{\begin{eqnarray}}
\newcommand{\ben}{\begin{eqnarray*}}
\newcommand{\en}{\end{eqnarray}}
\newcommand{\enn}{\end{eqnarray*}}

\newtheorem{theorem}{Theorem}[section]
\newtheorem{lemma}[theorem]{Lemma}

\newtheorem{remark}[theorem]{Remark}

\definecolor{rot}{rgb}{1,0,0}
\definecolor{hw}{rgb}{0,0,1}

\begin{document}
\renewcommand{\theequation}{\arabic{section}.\arabic{equation}}
\title{\bf
 Inverse source problem for a hyperbolic equation
 by Carleman estimates 
}
\author{ Suliang Si\thanks{School of Mathematics and Statistics, Shandong University of Technology,
Zibo, 255000, China ({\tt sisuliang@amss.ac.cn})}}
\date{}




\maketitle

\begin{abstract}
 In this article, we provide a modified argument for proving the conditional stability
of inverse source problem for a hyperbolic equation. Our method does not require any extension of solution with respect to time and
 therefore simplifies the existing proofs, which is widely applicable to various evolution equations.
\end{abstract}

%


\section{Introduction}
\noindent 

Let \( \Omega \subset \mathbb{R}^n \), \(n\geq 2\) be a bounded open domain with smooth boundary \( \partial \Omega \).
Denote by \( \nu = \nu(x) \) the unit outward normal vector to \( \partial \Omega \) at the point \( x \) and set \( \partial_\nu u = \nabla u \cdot \nu \).

We consider an inverse source problem for wave equation:
\begin{equation}\label{u}
\begin{cases}
\partial_t^2 u - \Delta u + \sum_{j=1}^n b_j(x)\partial_{x_j} u +d(x)\partial_tu+ c(x)u = R(x,t)f(x), & x \in \Omega, \quad 0 < t < T, \\
u(x, 0) = \partial_t u(x, 0) = 0, & x\in\Omega, \\
u = 0,  & (x,t)\in\partial \Omega \times (0,T),
\end{cases}
\end{equation}
where
\begin{equation}\label{b}
b_j, d, c \in L^\infty(\Omega), \quad \forall \, j \in \{1,\ldots,n\}.
\end{equation}
$R(x,t)f(x)$ accounts for the external force. 
For arbitrarily fixed \( x_0 \notin \overline{\Omega} \), we set
\begin{equation}\label{x}
\Gamma_0 := \{ x \in \partial \Omega|\, \, (x - x_0) \cdot \nu(x) \geq 0 \}.
\end{equation}
In this article, we consider the following inverse problem:
 Given $R(x,t)$,  determine the
 source $f(x)$, $x\in\Omega$ in (\ref{u}) by observation data $\partial_\nu u$ on a subboundary $\Sigma_0=\Gamma_0\times(0,T)$

We now state the main result for the inverse source problem.
\begin{theorem}\label{thm}
Assume that (\ref{b}) and (\ref{x}) hold, \(f\in L^2(\Omega)\), \( R \in H^1(0, T; L^\infty(\Omega)) \) and there exists a constant \( r_0 > 0 \) such that
\begin{equation}\label{r0}
|R(x, 0)| \geqslant r_0, \quad x \in \overline{\Omega}
\end{equation}
and
\begin{equation}\label{T0}
T > \left( \max_{x \in \overline{\Omega}} |x - x_0|^2 - \min_{x \in \overline{\Omega}} |x - x_0|^2 \right)^{\frac{1}{2}}.
\end{equation}
Then there exists a constant \( C > 0 \) such that

\[
\| f \|_{L^2(\Omega)} \leqslant C \| \partial_t \partial_\nu u \|_{L^2(\Sigma_0)}
\]
for each function \( u \).
\end{theorem}
\begin{remark}
The inverse source problem is not only interesting itself, but also is related to

\textbf{Coefficient inverse problem.} Determine $c(x)$ in 
\begin{equation}\label{w}
\begin{cases}
\partial_t^2 w(x,t) - \Delta w(x,t) + \sum_{j=1}^n b_j(x)\partial_{x_j} w +d(x)\partial_tw+ c(x)w = 0, & (x,t) \in Q, \\
w(x,0)=w_0, \quad \partial_tw(x,0)=w_1, & x\in\Omega,\\
w(x,t) = h, & (x,t) \in \partial\Omega \times (0,T).
\end{cases}
\end{equation}
where \( w_0 \), \( w_1 \), \(b_j\), \(d\) and \( h \) are given suitably. In fact, let \( w_j \) satisfy the above system with \( c_j \), \( j = 1, 2 \). Setting \( R := w_2 \), \( c := c_1 \), \( f := c_2 - c_1 \) and \( u := w_1 - w_2 \), we have the system (\ref{u}), and thus reduce the coefficient inverse problem to the inverse source problem for (\ref{u}).
\end{remark}

 The determination of coefficients for hyperbolic equations from boundary measurements has attracted much attention in recent years.
There is a considerable number of papers dealing with the uniqueness and stability in an
 inverse problem to identify unknown coefficients \cite{Bellassoued2017,JLY2017,K2002,Y1995,Y1999}.
In 1981 Bukhgeim and Klibanov
 \cite{BK1981} proposed a remarkable method based on a Carleman estimate and established the uniqueness for inverse
 problems of determining spatially varying coefficients for scalar partial differential equations. 
 Since then, many authors have successfully applied Carleman estimates to obtain Lipschitz stability or Hölder stability for hyperbolic and parabolic
 Cauchy problem with data on the lateral boundary as well as the inverse problems \cite{HIY,Beilina2012,Cannarsa2019a,Fu2019,Gölgeleyen2016,
 Imanuvilov2001,Imanuvilov2001a,Imanuvilov2003,Isakov1990,JLY2017}. 
In \cite{Imanuvilov1998} Imanuvilov and Yamamoto give a rigorous treatment for the uniqueness and stability
 by the Dirichlet data and the Neumann data on a sufficiently large part of the boundary 
over a sufficiently long time interval and proved a global and both-sided Lipschitz stability
 estimate.  
 estimate. 
 As for applications of Carleman estimates to inverse problems, we can refer to.  Most of these papers treat the determination of the coefficient \( p(x) \) in the zeroth-order term of a hyperbolic equation
\[
\partial_{t}^{2} u(t, x) - \Delta u(t, x) + p(x)u(t, x) = 0.
\]
 Especially, Huang,  Imanuvilov and Yamamoto  proposed an argument without the cut-off procedure for proving the stability for the inverse problems \cite{HIY}.

In the existing works, whenever one applied Carleman estimates, one needed to extend the solution about $t$, but this extension process does not yield any additional information. It is merely necessary for applying the Carleman estimate as it requires  $u(x,t)=\partial_tu(x,t)=0$, $t=\pm T$.
 In this article, we propose an argument without the extension of solution for proving the stability of the inverse source problem. The key is that we establish a new Carleman estimate, which is the main innovation of this paper.  

 Thearticle is composed of five sections. In Section \ref{02}, we present the key Carleman estimates. Sections \ref{pro} are devoted to the proofs of theorems \ref{thm}. In the appendix, we provide a detailed proof of the Carleman estimate (\ref{Carle}).

\section{Carleman estimates for the hyperbolic equation}\label{02}
 Throughout this article, $C, C_1>0$ denote generic constants which are independent of
 parameter $s > 0$.
 Let $Q=\Omega\times(0,T)$.

We first consider the following hyperbolic equation:
\begin{equation}\label{Car}
\begin{cases}
\partial_t^2 v(x,t) - \Delta v(x,t) + \sum_{j=1}^n b_j(x)\partial_{x_j} v(x,t) +d(x)v(x,t)+ c(x)v(x,t) = F(x,t), & (x,t) \in Q, \\
v(x,0)=0, \quad \partial_tv(x,0)=v_1, & x\in\Omega,\\
v(x,t) = 0, & (x,t) \in \partial\Omega \times (0,T).
\end{cases}
\end{equation}
For arbitrarily fixed \( x_0 \notin \overline{\Omega} \), \( \lambda >0 \), \( 0 \leqslant t_0 < T \) and \( 0 < \beta <1 \), we set
\[
\varphi(x,t) = \exp[{\lambda \psi(x,t)}], \quad \psi(x,t) = |x - x_0|^2 - \beta(t - t_0)^2, \quad (x,t) \in Q.
\tag{2.2}
\]

\begin{lemma}\label{lem}
Let the coefficients \( b_j, d, c \) satisfy (\ref{b}), \( \lambda >0 \) be sufficiently large and let \( \Gamma_0 \subset \partial \Omega \) be defined by (\ref{x}). Then there exist constants \( s_0 >0 \) and \( C >0 \) such that

\begin{equation}\label{Carle}
\begin{split}
 \int_{Q} \left( s|\nabla v|^2+s|\partial_t v|^2 + s^3 |v|^2 \right) e^{2s\varphi} dxdt &\leqslant C \left( \int_{Q} |F|^2 e^{2s\varphi} dxdt + \int_{\Gamma_0 \times (0,T)} s|\partial_\nu v|^2 e^{2s\varphi} \, dsdt \right. 
\\
& \left. + \int_{\Omega} \left( s|\nabla v(x,T)|^2+s|\partial_t v(x,T)|^2 + s^3 |v(x,T)|^2 \right) e^{2s\varphi(x,T)} \, dx \right). \\
\end{split}
\end{equation}
for all \( s > s_0 \) and \( v \in H^2(Q)\)  satisfying (\ref{Car}).
\end{lemma}
This is our most important result. Compared to the classical Carleman estimate, for example, to
 theorem 4.2 in \cite{Bellassoued2017}, we don't need the condition $\partial_tu(x,0)=0$ or to extend the solution in time.. We place the proof in the appendix.

\section{Proof of theorem 1.1}\label{pro}
In this section, we divide the proof of Theorem \ref{thm} into three steps.

\textbf{First step.}

By (\ref{u}) and \( R \in H^1(0,T;L^\infty(\Omega))\), setting \( y = \partial_t u \), we have

\begin{equation}\label{y}
\begin{cases}
\partial_t^2 y - \Delta y + \sum\limits_{j=1}^n b_j(x) \partial_{x_j} y +d(x)y +c(x) y = \partial_t R(x,t) f(x) \in L^2(Q), & x \in \Omega, \quad 0 < t < T, \\
y(x, 0) = 0, \quad \partial_t y(x, 0) = R(x, 0) f, & x\in\Omega, \\
y(x,t) = 0, &(x,t)\in \partial \Omega \times (0,T).
\end{cases}
\end{equation}
Set
\[
d_0 := \min_{x \in \overline{\Omega}} |x - x_0|, \quad d_1 := \max_{x \in \overline{\Omega}} |x - x_0|.
\]
Let \( t_0 = 0 \). Since (\ref{T0}) implies that \( T > \sqrt{d_1^2 - d_0^2} \), we choose \( \beta \in (0,1) \) sufficiently close to 1, such that
\begin{equation}\label{T}
T > \frac{\sqrt{d_1^2 - d_0^2}}{\sqrt{\beta}}.
\end{equation}
Therefore, we apply Lemma \ref{lem} to \( y \) in \( Q \):
\begin{equation}\label{yc}
\begin{aligned}
& \int_Q (s|\nabla y|^2 +s|\partial_t y|^2+ s^3|y|^2)\mathrm{e}^{2s\varphi}\mathrm{d}x\mathrm{d}t \\
&\leqslant C\left( J + \mathrm{e}^{C s}\int_{\Gamma \times (0, T)} |\partial_t \partial_\nu u|^2 \mathrm{d}\Sigma + s^3 \int_\Omega (|\nabla y(x, T)|^2 +|\partial_t y(x, T)|^2+ |y(x, T)|^2)\mathrm{e}^{2s\varphi(x, T)}\mathrm{d}x \right)
\end{aligned}
\end{equation}
for all \( s \geqslant s_0 \). Here we set
\[
J := \int_Q |\partial_t R|^2 |f|^2 \mathrm{e}^{2s\varphi}\mathrm{d}x\mathrm{d}t.
\]

\textbf{Second step.}

There exists a constant \( C > 0 \) such that
\begin{equation}\label{y0}
\int_\Omega |\partial_t y(x, 0)|^2 \mathrm{e}^{2s\varphi(x,0)} \mathrm{d}x\\
\leqslant C J + C \int_Q (s |\nabla_{x,t} y|^2 + |y|^2) \mathrm{e}^{2s\varphi} \mathrm{d}x\mathrm{d}t + \int_\Omega |\nabla_{x,t} y(x,T)|^2 \mathrm{e}^{2s\varphi(x,T)} \mathrm{d}x
\end{equation}
for all large \( s \geqslant 1 \),
which can be found in \cite{HIY}.

\textbf{Third step.}

We will complete the proof of theorem \ref{thm} by (\ref{yc}) and (\ref{y0}). The second equation in (\ref{y}) implies
\[
\partial_t y(\cdot, 0) = R(\cdot, 0) f \quad \text{in } \Omega.
\]
Therefore, by (\ref{r0}), estimate (\ref{y0}) yields
\begin{equation}\label{of}
\int_\Omega |f(x)|^2 \mathrm{e}^{2s\varphi(x,0)} \mathrm{d}x
\\
\leqslant C \left( J + \int_Q (s |\nabla_{x,t} y|^2 + |y|^2) \mathrm{e}^{2s\varphi} \mathrm{d}x\mathrm{d}t + \int_\Omega |\nabla_{x,t} y(x,T)|^2 \mathrm{e}^{2s\varphi(x,T)} \mathrm{d}x \right).
\end{equation}
We apply (\ref{yc}) to estimate the second term on the right-hand side and obtain
\begin{equation}\label{f}
\int_\Omega |f(x)|^2 \mathrm{e}^{2s\varphi(x,0)} \mathrm{d}x \leqslant C \left( J + \mathrm{e}^{C s} \|\partial_t \partial_\nu u\|^2_{L^2(\Gamma \times (0,T))} \right.
\\
\left. + \, s^3 \int_\Omega (|\nabla_{x,t} y(x,T)|^2 + |y(x,T)|^2) \mathrm{e}^{2s\varphi(x,T)} \mathrm{d}x \right)
\end{equation}
for sufficiently large \( s > 0 \).

On the other hand, there exists a constant \( s_1 > 1 \) such that
\[
J = o(1) \int_\Omega |f(x)|^2 \mathrm{e}^{2s\varphi(x,0)} \mathrm{d}x, \quad \forall \, s \geqslant s_1.
\]
Here \( o(1) \) denotes quantities which converge to 0 as \( s \to \infty \).

Therefore, choosing $s>0$ sufficiently large, we absorb the first term on the right-hand side
 of (\ref{f}) into the left-hand side:
\begin{equation}\label{f1}
\int_\Omega |f(x)|^2 \mathrm{e}^{2s\varphi(x,0)} \mathrm{d}x \leqslant C \left(\mathrm{e}^{C s} \|\partial_t \partial_\nu u\|^2_{L^2(\Gamma \times (0,T))} \right.
\\
\left. + \, s^3 \int_\Omega (|\nabla_{x,t} y(x,T)|^2 + |y(x,T)|^2) \mathrm{e}^{2s\varphi(x,T)} \mathrm{d}x \right)
\end{equation}
for sufficiently large \( s \). Here we apply the classical a priori estimate and we see
\[
\int_\Omega |\nabla_{x,t} y(x,T)|^2 \mathrm{d}x \leqslant C \|f\|^2_{L^2(\Omega)}.
\]
Moreover, the Poincaré inequality and the zero Dirichlet boundary condition for the function \( y(\cdot, T) \), yield
\[
\int_\Omega |y(x,T)|^2 \mathrm{d}x \leqslant C \int_\Omega |\nabla y(x,T)|^2 \mathrm{d}x.
\]
Hence
\begin{equation}\label{f00}
\begin{split}
\int_\Omega &(|\nabla_{x,t} y(x,T)|^2 + |y(x,T)|^2) \mathrm{e}^{2s\varphi(x,T)} \mathrm{d}x
\\
&\leqslant C \mathrm{e}^{2s \mathrm{e}^{\lambda (d_1^2 - \beta T^2)}} \int_\Omega (|\nabla_{x,t} y(x,T)|^2 + |y(x,T)|^2) \mathrm{d}x \leqslant C \mathrm{e}^{2s \mathrm{e}^{\lambda (d_1^2 - \beta T^2)}} \|f\|^2_{L^2(\Omega)}.
\end{split}
\end{equation}
On the other hand, we have
\[
\int_\Omega |f(x)|^2 \mathrm{e}^{2s\varphi(x,0)} \mathrm{d}x = \int_\Omega \mathrm{e}^{2s \mathrm{e}^{\lambda |x - x_0|^2}} |f(x)|^2 \mathrm{d}x \geqslant \mathrm{e}^{2s \mathrm{e}^{\lambda d_0^2}} \|f\|^2_{L^2(\Omega)}.
\]
Consequently (\ref{f1}) yields
\[
\|f\|^2_{L^2(\Omega)} \leqslant C s^3 \mathrm{e}^{-c_0 s} \|f\|^2_{L^2(\Omega)} + C \mathrm{e}^{C s} \|\partial_t \partial_\nu u\|^2_{L^2(\Gamma \times (0,T))},
\]
where \( c_0 = 2 \left( \mathrm{e}^{\lambda d_0^2} - \mathrm{e}^{\lambda d_1^2 - \lambda \beta T^2} \right) \). The inequality (\ref{T}) yields \( c_0 > 0 \). Finally, by noting \( \lim_{s \to \infty} s^3 \mathrm{e}^{-c_0 s} = 0 \), we absorb the first term on the right-hand side by taking sufficiently large \( s \). This proves Theorem \ref{thm}.

\section{Conclusion}
We have presented a stability result for the inverse source problem of a hyperbolic equation by Carleman estimates. Our approach improves Bukhgeim and Klibanov's proof by eliminating the cut-off procedure and time extension of the solution. A key point is that our Carleman estimates don't require $\partial v(x,0)=0$.  A possible continuation of this work is to investigate the fixed angle inverse scattering problems.

\section*{Acknowledgment}
The work of Suliang Si is supported by  the Shandong Provincial Natural Science Foundation (No. ZR2022QA111).

\section{Appendix. Carleman estimates}
Now, we provide a complete proof of Lemma \ref{lem}.
Write
\[
v'=\partial_t v, \quad v''=\partial_t^2 v,
\]
if there is no danger of confusion. Let $\Sigma=\partial\Omega\times(0,T)$.

It is sufficient to prove Lemma \ref{lem} in the case where \( b_j,d, c\equiv 0, \, \forall \, j \in \{1,\ldots,n\}  \). Indeed,  we assume that we already established the inequality
\begin{equation}\label{car000}
\begin{split}
\int_{Q} \left( s|\nabla v|^2 + s|\partial_t v|^2+ s^3 |v|^2 \right) e^{2s\varphi} dxdt &\leqslant C \left( \int_{Q} |F|^2 e^{2s\varphi} dxdt + \int_{\Gamma_0 \times (0,T)} s|\partial_\nu v|^2 e^{2s\varphi}dsdt \right. 
\\
& \left. + \int_{\Omega} \left( s|\nabla v(x,T)|^2 + s|\partial_t v(x,T)|^2+ s^3 |v(x,T)|^2 \right) e^{2s\varphi(x,T)} \, dx \right). \\
\end{split}
\end{equation}
Since the coefficients of \( b_j,d, c \) are in \( L^\infty(Q) \), we have
\begin{equation}
\begin{split}
|(\partial_t^2 - \Delta) v(x, t)|^2 &= |(Pv - P_1 v)(x, t)|^2 \leq 2|Pv(x, t)|^2 + 2|P_1 v(x, t)|^2\\
&\leq 2|Pv(x, t)|^2 + C \left( |v(x, t)|^2 + |\nabla v(x, t)|^2  + |v'(x, t)|^2 \right), \quad (x, t) \in Q,
\end{split}
\end{equation}
where \[Pv=\partial_t^2 v - \Delta v \] and
\[P_1v=- \sum_{j=1}^n b_j(x)\partial_{x_j} v--d(x)\partial_tv- c(x)v.\]
By choosing $s$ large, one can absorb the term
\[
\int_Q \left( |v(x, t)|^2 + |\nabla v(x, t)|^2 + |v'(x, t)|^2 \right) e^{2s\varphi} dx dt
\]
into the left-hand side of the Carleman estimate (\ref{car000}).
Thus we only need to prove Lemma \ref{lem} when \( b_j,d, c\equiv 0, \, \forall \, j \in \{1,\ldots,n\}  \).

Let 
\[
z(x, t) = e^{s \varphi} v(x, t), \quad (x, t) \in Q.
\]
Define 
\[
P_s^+ z = z'' - \Delta_{\mathrm{g}} z + s^2 \left( | \varphi' |^2 - | \nabla_{\mathrm{g}} \varphi |^2 \right) z
\]
and
\[
P_s^- z = -2s \left( z' \varphi' - \langle \nabla_{\mathrm{g}} z, \nabla_{\mathrm{g}} \varphi \rangle \right) - s \left( \varphi'' - \Delta_{\mathrm{g}} \varphi \right) z.
\]
Then
\begin{equation}
e^{s\varphi}(\partial_t^2-\Delta)v=e^{s\varphi}(\partial_t^2-\Delta)(e^{-s \varphi}z )=P_s^+ z+P_s^- z=e^{s\varphi}F.
\end{equation}
Therefore
\begin{equation}
P_sz:=P_s^+ z+P_s^- z=e^{s\varphi}F
\end{equation}
Let \( \|\cdot\| \) and \( (\cdot, \cdot) \) denote the norm \( \|\cdot\|_{L^2(Q)} \) and the scalar product \( (\cdot, \cdot)_{L^2(Q)} \) in \( L^2(Q) \).
With the above notations, we have
\[\|P_sz\|=\| P_s^+ z \|^2 + \| P_s^- z \|^2 + 2 \left( P_s^+ z, P_s^- z \right) = \| e^{s \varphi} F \|^2.
\]
Let us compute the term \( 2 \left( P_s^+ z, P_s^- z \right) \). For this, we will expand the six terms appearing in \( \left( P_s^+ z, P_s^- z \right) \) and integrate by parts repeatedly with respect to \( (x, t) \).

\begin{lemma}\label{Alem0}
   For any \(z \in H^2(Q)\) satisfying
    \[
        z(x,0)= 0, 
    \]
 the following identity holds true:
\begin{equation}\nonumber
\begin{split}
(P_s^+ z, P_s^- z) &= s\int_\Omega \big(\varphi' |z'|^2\big)(0)dx-s\int_\Omega \varphi'(T) |z'(T)|^2dx -s \int_\Omega \big(z z' \left( \varphi'' - \Delta \varphi \right)\big)(T)  dx
\\
&+\frac{s}{2} \int_\Omega \big(|z|^2 \left( \partial_t^2 - \Delta \right) \varphi'\big)(T) dx -s \int_\Omega \big( | \nabla z |^2 \varphi'\big)(T) dx -s^3 \int_\Omega \big(|z|^2  \varphi' \left( |\varphi'|^2 - |\nabla \varphi|^2 \right)\big)(T) dx
\\
&
+2s \int_Q \left( \varphi'' |z'|^2 - 2 z' \langle \nabla z, \nabla_g \varphi' \rangle +
+ \nabla^2 \varphi (\nabla z, \nabla z) \right) dx dt
\\
&+ 2s^3 \int_Q |z|^2 \left( |\varphi'|^2 \varphi'' + \nabla^2 \varphi (\nabla \varphi, \nabla \varphi) - 2 \varphi' \langle \nabla \varphi, \nabla \varphi' \rangle \right) dx dt
\\
&- \frac{s}{2} \int_Q |z|^2 \left( \partial_t^2 - \Delta \right)^2 \varphi dx dt + \mathcal{B}_0.
\end{split}
\end{equation}
where \( B_0 \) is given by
\begin{equation}\label{AB}
\begin{split}
\mathcal{B}_0 =& s \int_{\Sigma} \left[ \partial_\nu \varphi |\nabla z|^2 - 2 \langle \nabla z, \nabla \varphi \rangle \partial_\nu z + \left(2\varphi' z' \partial_\nu z - |z'|^2 \partial_\nu \varphi\right) \right]  dsdt\\
&+ s \int_{\Sigma} \left[ z \partial_\nu z (\varphi'' - \Delta \varphi) + s^2 \partial_\nu \varphi |z|^2 \left(|\varphi'|^2 - |\nabla \varphi|^2\right) - \frac{1}{2} |z|^2 \partial_\nu (\varphi'' - \Delta \varphi) \right]  dsdt    
\end{split}
\end{equation}   
and
\[
\nabla^{2}\varphi(\nabla \varphi, \nabla \varphi) = \sum_{i=1}^{n}\sum_{j=1}^{n}\frac{\partial^{2}\varphi}{\partial x_{i}\partial x_{j}} \cdot \frac{\partial \varphi}{\partial x_{i}} \cdot \frac{\partial \varphi}{\partial x_{j}}.
\]
\end{lemma}

\begin{proof}
 By the definitions of \( P_s^{-} \) and \( P_s^{+} \), we have
\begin{equation}\label{I6}
\begin{split}
\left( P_s^{+} z, P_s^{-} z \right)
&= -2s \int_Q z'' \left( z' \varphi' - \langle \nabla z, \nabla \varphi \rangle \right) dx dt - s \int_Q z'' \left( \varphi'' - \Delta \varphi \right) z dx dt
\\
&+ 2s \int_Q \Delta z \left( z' \varphi' - \langle \nabla z, \nabla \varphi \rangle \right) dx dt + s \int_Q \Delta z \left( \varphi'' - \Delta \varphi \right) z dx dt
\\
&- 2s^3 \int_Q \left( |\varphi'|^2 - |\nabla \varphi|^2 \right) z \left( z' \varphi' - \langle \nabla z, \nabla \varphi \rangle \right) dx dt
\\
&- s^3 \int_Q \left( |\varphi'|^2 - |\nabla \varphi|^2 \right) \left( \varphi'' - \Delta \varphi \right) |z|^2 dx dt := \sum_{k=1}^6 I_k.
\end{split}
\end{equation}
We calculate the six terms \( I_k \), \( k = 1, \ldots, 6 \) by integrating by parts with respect to \( (x, t) \).

First one easily sees that
\begin{equation}\nonumber
\begin{split}
I_1 &= -s \int_Q \varphi' \partial_t (|z'|^2) dx dt - s \int_Q \langle \nabla (|z'|^2), \nabla \psi \rangle dx dt
- 2s \int_Q z' \langle \nabla \varphi', \nabla z \rangle dx dt
\\
&= s\int_\Omega \big(\varphi' |z'|^2\big)(0)dx-s\int_\Omega \varphi'(T) |z'(T)|^2dx+s \int_Q |z'|^2 (\varphi'' + \Delta \varphi) dx dt 
\\
&- 2s \int_Q z' \langle \nabla_g \varphi', \nabla z \rangle dx dt
- s \int_\Sigma |z'|^2 \partial_\nu \varphi ds dt.
\end{split}
\end{equation}
Integration by parts yields
\begin{equation}\nonumber
\begin{split}
I_2 &= -s \int_Q z'' \left( \varphi'' - \Delta \varphi \right) z dx dt
\\
&= -s \int_\Omega \big(z z' \left( \varphi'' - \Delta \varphi \right)\big)(T)  dx +s \int_Q \left( \varphi'' - \Delta \varphi \right) |z'|^2 dx dt + \frac{s}{2} \int_Q \partial_t (|z|^2) \left( \partial_t^2 - \Delta \right) \varphi' dx dt
\\
&=-s \int_\Omega \big(z z' \left( \varphi'' - \Delta \varphi \right)\big)(T)  dx+\frac{s}{2} \int_\Omega \big(|z|^2 \left( \partial_t^2 - \Delta \right) \varphi'\big)(T) dx \\
&\quad+ s \int_Q \left( \varphi'' - \Delta_g \varphi \right) |z|^2 dx dt - \frac{s}{2} \int_Q |z|^2 \left( \partial_t^2 - \Delta \right) \varphi'' dx dt.
\end{split}
\end{equation}
Furthermore, by Green's formula and integration by parts, we obtain
\begin{equation}\nonumber
\begin{split}
I_3 &= 2s \int_Q \Delta z \left( z' \varphi' - \langle \nabla z, \nabla \varphi \rangle \right) dx dt
\\
&= -2s \int_Q \left( \langle \nabla z, \nabla z' \rangle \varphi' + \langle \nabla z, \nabla \varphi' \rangle z' - \langle \nabla z, \nabla \left( \langle \nabla z, \nabla \varphi \rangle \right) \right) dx dt
\\
&+ 2s \int_\Sigma \partial_\nu z \left( z' \varphi' - \langle \nabla_g z, \nabla_g \varphi \rangle \right) ds dt
\\
&=-s \int_\Omega \big( | \nabla z |^2 \varphi'\big)(T) dx + s \int_Q \left( | \nabla z |^2 \varphi'' - 2 \langle \nabla z, \nabla \varphi' \rangle z' + 2 \langle \nabla z, \nabla \left( \langle \nabla z, \nabla \varphi \rangle \right) \right) dx dt
\\
&+ 2s \int_\Sigma \partial_\nu z \left( z' \varphi' - \langle \nabla z, \nabla \varphi \rangle \right) ds dt.
\end{split}
\end{equation}
Applying
\[
\langle \nabla z, \nabla_g \left( \langle \nabla z, \nabla \varphi \rangle \right) \rangle = \nabla^2 \varphi \left( \nabla z, \nabla z \right) + \frac{1}{2} \langle \nabla \varphi, \nabla (| \nabla z |^2 ),
\]
we conclude that
\begin{equation}\nonumber
\begin{split}
I_3 &= -s \int_\Omega \big( | \nabla z |^2 \varphi'\big)(T) dx+ s \int_Q \left( |\nabla z|^2 \left( \varphi'' - \Delta \varphi \right) - 2 z' \langle \nabla z, \nabla \varphi' \rangle + 2 \nabla_g^2 \varphi (\nabla z, \nabla z) \right) dx dt
\\
&
+ s \int_\Sigma \left( 2 \partial_\nu z \left( z' \varphi' - \langle \nabla z, \nabla \varphi \rangle \right) + \partial_\nu \varphi |\nabla z|^2 \right) ds dt.
\end{split}
\end{equation}
On the other hand,
\begin{equation}\nonumber
\begin{split}
I_4 &= s \int_Q \Delta z \left( \varphi'' - \Delta \varphi \right) z dx dt
\\
&= -s \int_Q \left( |\nabla z|^2 \left( \varphi'' - \Delta \varphi \right) + \frac{1}{2} \langle \nabla (|z|^2), \nabla \left( \varphi'' - \Delta \varphi \right) \rangle \right) dx dt
+ s \int_\Sigma \partial_\nu z \left( \varphi'' - \Delta \varphi \right) z ds dt
\\
&= -s \int_Q \left( |\nabla z|^2 \left( \varphi'' - \Delta \varphi \right) - \frac{1}{2} |z|^2 \Delta \left( \varphi'' - \Delta \varphi \right) \right) dx dt.
\\
&+ s \int_\Sigma \left( \partial_\nu z \left( \varphi'' - \Delta \varphi \right) z - \frac{1}{2} |z|^2 \partial_\nu ( \varphi'' - \Delta \varphi ) \right) ds dt.
\end{split}
\end{equation}
Next, we have
\begin{equation}\nonumber
\begin{split}
I_5 &= -2s^3 \int_Q \left( |\varphi'|^2 - |\nabla \varphi|^2 \right) z \left( z' \varphi' - \langle \nabla \varphi, \nabla z \rangle \right) dx dt
\\
&= -s^3 \int_Q \left( \partial_t (|z|^2) \varphi' - \langle \nabla (|z|^2), \nabla \varphi \rangle \right) \left( |\varphi'|^2 - |\nabla \varphi|^2 \right) dx dt
\\
&= -s^3 \int_\Omega \big(|z|^2  \varphi' \left( |\varphi'|^2 - |\nabla \varphi|^2 \right)\big)(T) dx + s^3 \int_Q |z|^2 \left( \varphi'' - \Delta \varphi \right) \left( |\varphi'|^2 - |\nabla \varphi|^2 \right) dx dt
\\
&+ s^3 \int_Q |z|^2 \left( \varphi' \partial_t (|\varphi'|^2 - |\nabla \varphi|^2) - \langle \nabla \varphi, \nabla (|\varphi'|^2 - |\nabla \varphi|^2) \rangle \right) dx dt
\\
&+ s^3 \int_\Sigma \partial_\nu \varphi |z|^2 \left( |\varphi'|^2 - |\nabla \varphi|^2 \right) ds dt
\\
&=-s^3 \int_\Omega \big(|z|^2  \varphi' \left( |\varphi'|^2 - |\nabla \varphi|^2 \right)\big)(T) dx+ s^3 \int_Q |z|^2 \left( \varphi'' - \Delta \varphi \right) \left( |\varphi'|^2 - |\nabla \varphi|^2 \right) dx dt
\\
&+ 2s^3 \int_Q |z|^2 \left( |\varphi'|^2 \varphi'' + \nabla^2 \varphi (\nabla \varphi, \nabla \varphi) - 2 \varphi' \langle \nabla \varphi, \nabla \varphi' \rangle \right) dx dt
\\
&+ s^3 \int_\Sigma \partial_\nu \varphi |z|^2 \left( |\varphi'|^2 - |\nabla \varphi|^2 \right) ds dt.
\end{split}
\end{equation}
Finally,
\[
I_6 = -s^3 \int_Q |z|^2 \left( \varphi'' - \Delta_g \varphi \right) \left( |\varphi'|^2 - |\nabla_g \varphi|^2 \right) dx dt.
\]
Then by (\ref{I6}), we obtain
\begin{equation}\label{IB}
\begin{split}
(P_s^+ z, P_s^- z) &= s\int_\Omega \big(\varphi' |z'|^2\big)(0)dx-s\int_\Omega \varphi'(T) |z'(T)|^2dx -s \int_\Omega \big(z z' \left( \varphi'' - \Delta \varphi \right)\big)(T)  dx
\\
&+\frac{s}{2} \int_\Omega \big(|z|^2 \left( \partial_t^2 - \Delta \right) \varphi'\big)(T) dx -s \int_\Omega \big( | \nabla z |^2 \varphi'\big)(T) dx -s^3 \int_\Omega \big(|z|^2  \varphi' \left( |\varphi'|^2 - |\nabla \varphi|^2 \right)\big)(T) dx
\\
&
+2s \int_Q \left( \varphi'' |z'|^2 - 2 z' \langle \nabla z, \nabla \varphi' \rangle +
+ \nabla^2 \varphi (\nabla z, \nabla z) \right) dx dt
\\
&+ 2s^3 \int_Q |z|^2 \left( |\varphi'|^2 \varphi'' + \nabla^2 \varphi (\nabla \varphi, \nabla \varphi) - 2 \varphi' \langle \nabla \varphi, \nabla \varphi' \rangle \right) dx dt
\\
&- \frac{s}{2} \int_Q |z|^2 \left( \partial_t^2 - \Delta \right)^2 \varphi dx dt + \mathcal{B}_0.
\end{split}
\end{equation}
where $B_0$ is given by (\ref{AB}). Thus the proof of Lemma \ref{Alem0} is complete.

\textbf{Second Step.} 

In this step, we want to prove a lower bound of \( (P_s^+ z, P_s^- z) \). To this end we decompose the right-hand side of (\ref{IB}) as
\[
(P_s^+ z, P_s^- z) = J_1 + J_2 + J_3 +J_4+ \mathcal{B}_0, 
\]

where
\begin{equation}\label{J1}
J_1 = 2s \int_Q \left( \varphi'' |z'|^2 - 2 z' \langle \nabla z, \nabla \varphi' \rangle + \nabla^2 \varphi (\nabla z, \nabla z) \right) dx dt,
\end{equation}
\[
J_2 = 2s^3 \int_Q |z|^2 \left( |\varphi'|^2 \varphi'' + \nabla^2 \varphi (\nabla \varphi, \nabla \varphi) - 2 \varphi' \langle \nabla \varphi, \nabla \varphi' \rangle \right) dx dt,
\]
\begin{equation}\label{J3}
J_3 = -\frac{s}{2} \int_Q |z|^2 \left( \partial_t^2 - \Delta \right)^2 \varphi dx dt,
\end{equation}
\begin{equation}\nonumber
\begin{split}
J_4&=s\int_\Omega \big(\varphi' |z'|^2\big)(0)dx-s\int_\Omega \varphi'(T) |z'(T)|^2dx -s \int_\Omega \big(z z' \left( \varphi'' - \Delta \varphi \right)\big)(T)  dx
\\
&+\frac{s}{2} \int_\Omega \big(|z|^2 \left( \partial_t^2 - \Delta \right) \varphi'\big)(T) dx -s \int_\Omega \big( | \nabla z |^2 \varphi'\big)(T) dx -s^3 \int_\Omega \big(|z|^2  \varphi' \left( |\varphi'|^2 - |\nabla \varphi|^2 \right)\big)(T) dx.
\end{split}
\end{equation}
Set
\[
b(\psi)(x, t) = |\psi'(x, t)|^2 - |\nabla_g \psi(x, t)|^2.
\]

\begin{lemma}\label{Alem1}
There exists a constant \(C > 0\) such that for any \(\varepsilon > 0\), there exists a constant \(C_{\varepsilon} > 0\) such that
\begin{equation}\nonumber
\begin{split}
J_1 + 2(1 + \beta)\mathcal{B}_1 +&2(1 +\beta) \int_\Omega \big(s\lambda\varphi z'\big)(T) dx -(1 +\beta) \int_\Omega \big((s\lambda\varphi)' |z|^2\big)(T) dx
\\
& \geq 2(1 - \beta) \int_Q s\lambda\varphi \left( |z'|^2 + |\nabla z|^2 \right) dx dt
\\
&
- C \left( \int_Q (s\lambda\varphi)^3 |b(\psi)| |z|^2 dx dt + C_{\varepsilon} \int_Q (s\lambda\varphi)^2 |z|^2 dx dt + \varepsilon \left\| P_s^+ z \right\|^2 \right),
\end{split}
\end{equation}
where \(\mathcal{B}_1\) is given by
\[
\mathcal{B}_1 = \int_{\Sigma} s\lambda\varphi \left( \lambda \partial_{\nu} \psi |z|^2 - z \partial_{\nu} z \right) ds dt.
\]
\end{lemma}
\begin{proof}
By expanding (\ref{J1}), we have
\begin{equation}\nonumber
\begin{split}
J_1 &= 2s \int_Q \lambda \varphi ( (\psi'' + \lambda |\psi'|^2) |z'|^2 - 2 \lambda \psi' z' \langle \nabla z, \nabla \psi \rangle 
\\
&\quad+ \nabla^2 \psi (\nabla z, \nabla z) + \lambda | \langle \nabla z, \nabla \psi \rangle |^2 ) dx dt
\\
&= 2 \int_Q s\lambda\varphi \left( \psi'' |z'|^2 + \nabla^2 \psi (\nabla z, \nabla z) + \lambda \left( \psi' z' - \langle \nabla z, \nabla \psi \rangle \right)^2 \right) dx dt.
\end{split}
\end{equation}
Then
\begin{equation}\label{J1}
\begin{split}
J_1 &\geq 2 \int_Q s\lambda\varphi \left( \nabla^2 \psi (\nabla z, \nabla z) + \psi'' |z'|^2 \right) dx dt
\\
& \geq 4 \int_Q s\lambda\varphi |\nabla z|^2 dx dt - 4 \beta \int_Q s\lambda\varphi |z'|^2 dx dt.
\end{split}
\end{equation}
Next, multiplying the hyperbolic equation in \( z \) by \( \sigma z \) and integrating by parts, we have
\begin{equation}\nonumber
\begin{split}
\int_Q P_s^+ z (s\lambda\varphi z) dx dt &=\int_\Omega \big(s\lambda\varphi z'\big)(T) dx  -\int_Q s\lambda\varphi |z'|^2 dx dt + \int_Q s\lambda\varphi |\nabla_g z|^2 dx dt - \frac{1}{2} \int_Q (s\lambda\varphi)' (\partial_t |z|^2) dx dt
\\
&+ \frac{1}{2} \int_Q \langle \nabla s\lambda\varphi, \nabla (|z|^2) \rangle dx dt + \int_Q (s\lambda\varphi)^3 b(\psi) |z|^2 dx dt - \int_\Sigma s\lambda\varphi z \partial_\nu z ds dt
\\
&= \int_\Omega \big(s\lambda\varphi z'\big)(T) dx -\frac{1}{2}\int_\Omega \big((s\lambda\varphi)' |z|^2\big)(T) dx
\\
& -\int_Q s\lambda\varphi |z'|^2 dx dt + \int_Q s\lambda\varphi |\nabla z|^2 dx dt + \frac{1}{2} \int_Q \big((s\lambda\varphi)'' - \Delta (s\lambda\varphi)\big) |z|^2 dx dt
\\
&+ \int_Q (s\lambda\varphi)^3 b(\psi) |z|^2 dx dt + \underbrace{\int_\Sigma s\lambda\varphi \left( \gamma \partial_\nu \psi |z|^2 - z \partial_\nu z \right) ds dt}_{\mathcal{B}_1}.
\end{split}
\end{equation}
Since
\[
(s\lambda\varphi)'' - \Delta (s\lambda\varphi) = s\lambda^2\varphi\left( \psi'' - \Delta \psi \right) + s\lambda^3\varphi b(\psi),
\]
we deduce that for any \( \varepsilon > 0 \), there exists \( C_\varepsilon > 0 \) such that
\begin{equation}\label{J10}
\begin{split}
&\left| \int_Q s\lambda\varphi |z'|^2 dx dt - \mathcal{B}_1 -\int_\Omega \big(s\lambda\varphi z'\big)(T) dx +\frac{1}{2}\int_\Omega \big((s\lambda\varphi)' |z|^2\big)(T) dx\right| 
\\
&\leq \int_Q (s\lambda\varphi)^3 |z|^2 |b(\psi)| dx dt + \varepsilon \left\| P_s^+ z \right\|^2
\\
&+ \int_Q s\lambda\varphi|\nabla_g z|^2 dx dt + C_\varepsilon \int_Q (s\lambda\varphi)^2 |z|^2 dx dt. 
\end{split}
\end{equation}
Combining (\ref{J1}) and (\ref{J10}), we obtain
\begin{equation}\label{J12}
\begin{split}
J_1 + 4\beta \mathcal{B}_1&+4\beta\int_\Omega \big(s\lambda\varphi z'\big)(T) dx -2\beta\int_\Omega \big((s\lambda\varphi)' |z|^2\big)(T) dx  \\
&\geq 4(1 - \beta) \int_Q s\lambda\varphi |\nabla z|^2 dx dt
\\
&- C \left( \int_Q (s\lambda\varphi)^3 |z|^2 |b(\psi)| dx dt + \varepsilon \left\| P_s^+ z \right\|^2 + C_\varepsilon \int_Q (s\lambda\varphi)^2 |z|^2 dx dt \right).
\end{split}
\end{equation}
Using (\ref{J10}) again, we have
\begin{equation}\label{J13}
\begin{split}
2(1 - \beta) \int_Q s\lambda\varphi |z'|^2 dx dt - & C \left( \int_Q (s\lambda\varphi)^3 |b(\psi)| |z|^2 dx dt\right)
+ C_{\varepsilon} \int_Q (s\lambda\varphi)^2 |z|^2 dx dt + \varepsilon \left\| P_s^+ z \right\|^2 
\\
&\leq 2(1 - \beta) \int_Q s\lambda\varphi |\nabla z|^2 dx dt + 2(1 - \beta) \mathcal{B}_1
\\
&+2(1 - \beta) \int_\Omega \big(s\lambda\varphi z'\big)(T) dx -(1 - \beta) \int_\Omega \big((s\lambda\varphi)' |z|^2\big)(T) dx. 
\end{split}
\end{equation}
Combining (\ref{J13}) and (\ref{J12}), we complete the proof of Lemma \ref{Alem1}
\end{proof}

\begin{lemma} 
The following inequality holds
\[
J_2 \geq 2\lambda \int_Q (s\lambda\varphi)^3 (b(\psi))^2 |z|^2 dx dt + 4 \int_Q (s\lambda\varphi)^3 \left( |\nabla \psi|^2 - \beta |\psi'|^2 \right) |z|^2 dx dt.
\]
\end{lemma}
\begin{proof}
Expanding $\varphi'$ and $\varphi''$, we have
\begin{equation}\nonumber
\begin{split}
J_2 =& 2s^3 \int_Q (s\lambda^2\varphi)^3 |\psi'|^2 \left( \psi'' + \lambda |\psi'|^2 \right) |z|^2 dx dt
- 4s^3 \int_Q \lambda^4 \varphi^3 |\psi'|^2 |\nabla \psi|^2 |z|^2 dx dt
\\
&+ 2s^3 \int_Q \lambda^3 \varphi^3 \left( \nabla^2 \psi (\nabla \psi, \nabla \psi) + \gamma |\nabla \psi|^4 \right) |z|^2 dx dt
\\
&= 2 \int_Q (s\lambda\varphi)^3 \left( \psi'' |\psi'|^2 |z|^2 + \nabla^2 \psi (\nabla \psi, \nabla \psi) |z|^2 \right) dx dt
+ 2\lambda \int_Q (s\lambda\varphi)^3 (b(\psi))^2 |z|^2 dx dt
\\
&\geq 2\lambda \int_Q (s\lambda\varphi)^3 (b(\psi))^2 |z|^2 dx dt + 4 \int_Q (s\lambda\varphi)^3 \left(  |\nabla \psi|^2 - \beta |\psi'|^2 \right) |z|^2 dx dt.
\end{split}
\end{equation}
This completes the proof of the lemma.
\end{proof}
On the other hand, by (\ref{J3}), we obtain
\begin{equation}\label{J30}
|J_3| \leq C \lambda^2 \int_Q s\lambda\varphi|z|^2 dx dt \leq C \lambda \int_Q (s\lambda\varphi)^2 |z|^2 dx dt. \tag{4.20}
\end{equation}
Now (\ref{J30}) and Lemmas \ref{Alem1} and \ref{Alem0}, yield
\begin{lemma}
There exists a constant \(C > 0\) such that for any \(\varepsilon > 0\) there exists \(C_\varepsilon > 0\) such that
\begin{equation}\nonumber
\begin{split}
J_1 + J_2 + J_3 + 2(1 + \beta)\mathcal{B}_1+& 2(1 +\beta) \int_\Omega \big(s\lambda\varphi z'\big)(T) dx -(1 +\beta) \int_\Omega \big((s\lambda\varphi)' |z|^2\big)(T) dx
\\
 &\geq 2(1 - \beta) \int_Q s\lambda\varphi \left( |\nabla z|^2 + |z'|^2 \right) dx dt
\\
&+ 2\lambda \int_Q (s\lambda\varphi)^3 (b(\psi))^2 |z|^2 dx dt + 4 \int_Q (s\lambda\varphi)^3 \left( \varrho |\nabla \psi|^2 - \beta |\psi'|^2 \right) |z|^2 dx dt
\\
&- C \left( \int_Q (s\lambda\varphi)^3 |z|^2 |b(\psi)| dx dt + \varepsilon \left\| P_s^+ z \right\|^2 + C_\varepsilon \lambda \int_Q (s\lambda\varphi)^2 |z|^2 dx dt \right).
\end{split}
\end{equation}
\end{lemma}
\textbf{Third Step.}
Set
\begin{equation}\label{B}
\begin{split}
\mathcal{B} &= \int_{\Sigma} s\lambda\varphi \left( \partial_{\nu} \psi_0 |\nabla z|^2 - 2 \langle \nabla z, \nabla \psi_0 \rangle \partial_{\nu} z \right) ds dt
+ \int_{\Sigma} s\lambda\varphi \left( 2 \psi' z' \partial_{\nu} z - |z'|^2 \partial_{\nu} \psi_0 \right) ds dt
\\
&+ \int_{\Sigma} \big\{ s\lambda\varphi ( z \partial_{\nu} z \left( -2\beta - \Delta \psi + \gamma b(\psi) \right) + (s\lambda\varphi)^2 \partial_{\nu} \psi_0 |z|^2 b(\psi) )
+ \frac{1}{2} |z|^2 \partial_{\nu} ( \Delta \psi_0 + \lambda |\nabla \psi_0|^2 ) \big\} ds dt
\\
&- 2(1 + \beta) \int_{\Sigma} s\lambda\varphi \left( \lambda \partial_{\nu} \psi |z|^2 - z \partial_{\nu} z \right) ds dt. \end{split}
\end{equation}
Since \(\beta < 1\), for small \(\eta > 0\), we have
\[
\beta(1 + \eta) < 1.
\]
Denote
\[
Q^{\eta} = \left\{ (x, t) \in Q| \ |b(\psi)(x, t)| \leq \eta |\nabla \psi(x, t)|^2 \right\}.
\]
Then
\begin{equation}\label{Q}
\begin{split}
J_1 + J_2 + J_3 + 2(1 + \beta)\mathcal{B}_1 &\geq 2(1 - \beta) \int_Q s\lambda\varphi \left( |\nabla z|^2 + |z'|^2 \right) dx dt
\\
&+ 2\gamma \int_{Q \setminus Q^{\eta}} (s\lambda\varphi)^3 (b(\psi))^2 |z|^2 dx dt + 4(1 - \beta(1 + \eta)) \int_{Q^{\eta}} (s\lambda\varphi)^3 |z|^2 |\nabla \psi|^2 dx dt
\\
&- C \bigg( \eta \int_{Q^{\eta}} (s\lambda\varphi)^3 |z|^2 dx dt + \int_{Q \setminus Q^{\eta}} (s\lambda\varphi)^3 |z|^2 dx dt + \varepsilon \left\| P_s^+ z \right\|^2
\\
&+ C_{\varepsilon} \lambda \int_Q (s\lambda\varphi)^2 |z|^2 dx dt \bigg). 
\end{split}
\end{equation}
Using (\ref{Q}), we obtain
\begin{equation}\nonumber
\begin{split}
J_1 + J_2 + J_3 +& 2(\varrho + \beta)\mathcal{B}_1+2(1 +\beta) \int_\Omega \big(s\lambda\varphi z'\big)(T) dx -(1 +\beta) \int_\Omega \big((s\lambda\varphi)' |z|^2\big)(T) dx
\\
& \geq 2(1-\beta) \int_Q s\lambda\varphi \left( |\nabla z|^2 + |z'|^2 \right) dx dt
\\
&+ 2\gamma \eta^2 C_1 \int_{Q \setminus Q^{\eta}} (s\lambda\varphi)3 |z|^2 dx dt + C_2(1 - \beta(1 + \eta)) \int_{Q^{\eta}} (s\lambda\varphi)^3 |z|^2 dx dt
\\
&- C \left( \eta \int_{Q^{\eta}} (s\lambda\varphi)^3 |z|^2 dx dt + \int_{Q \setminus Q^{\eta}} (s\lambda\varphi)^3 |z|^2 dx dt + \varepsilon \left\| P_s^+ z \right\|^2 + \gamma \int_Q (s\lambda\varphi)^2 |z|^2 dx dt \right)
\\
&\geq \delta 2(1-\beta) \int_Q s\lambda\varphi \left( |\nabla z|^2 + |z'|^2 \right) dx dt + (2\gamma \eta^2 C_1 - C) \int_{Q \setminus Q^{\eta}} (s\lambda\varphi)^3 |z|^2 dx dt
\\
&+ (C_2(1 - \beta(1 + \eta)) - \eta C) \int_{Q^{\eta}} (s\lambda\varphi)^3 |z|^2 dx dt - C \left( \varepsilon \left\| P_s^+ z \right\|^2 + \lambda \int_Q (s\lambda\varphi)^2 |z|^2 dx dt \right).
\end{split}
\end{equation}
Then for small \(\eta\) and large \(\gamma \geq \gamma_0\) and \(s \geq s_0\) we obtain
\begin{equation}\nonumber
\begin{split}
J_1 + J_2 + J_3 + 2(1 + \beta)\mathcal{B}_1+&2(1 +\beta) \int_\Omega \big(s\lambda\varphi z'\big)(T) dx -(1 +\beta) \int_\Omega \big((s\lambda\varphi)' |z|^2\big)(T) dx
\\& \geq C \int_Q s\lambda\varphi\left( |\nabla z|^2 + |z'|^2 + \sigma^2 |z|^2 \right) dx dt - \frac{1}{4} \left\| P_s^+ z \right\|^2.
\end{split}
\end{equation}
By (\ref{IB}), we have
\begin{equation}\nonumber
\begin{split}
2 \left( P_s^+ z, P_s^- z \right)-2J_4 -2\mathcal{B}+&2(1 +\beta) \int_\Omega \big(s\lambda\varphi z'\big)(T) dx -(1 +\beta) \int_\Omega \big((s\lambda\varphi)' |z|^2\big)(T) dx- 2 \mathcal{B}
\\& \geq C \int_Q s\lambda\varphi \left( |\nabla z|^2 + |z'|^2 + \sigma^2 |z|^2 \right) dx dt - \frac{1}{2} \left\| P_s^+ z \right\|^2,
\end{split}
\end{equation}
where
\[
\mathcal{B} = \mathcal{B}_0 - 2(1 + \beta)\mathcal{B}_1.
\]
Then there exists a constant \( s_0 > 0 \) such that for any \( s \geq s_0 \) we have
\begin{equation}\label{Ps}
\begin{split}
\| P_s z \|^2 -2J_4- 2 \mathcal{B}+&2(1 +\beta) \int_\Omega \big(s\lambda\varphi z'\big)(T) dx -(1 +\beta) \int_\Omega \big((s\lambda\varphi)' |z|^2\big)(T) dx 
\\&\geq C \int_Q s\lambda\varphi \left( |\nabla z|^2 + |z'|^2 + (s\lambda\varphi)^2 |z|^2 \right) dx dt.
\end{split}
\end{equation}
Since \(\int_\Omega \big(\varphi' |z'|^2\big)(0)dx= 0\), we obtain
\begin{equation}\nonumber
\begin{split}
J_4&=-s\int_\Omega \varphi'(T) |z'(T)|^2dx -\frac{s}{2} \int_\Omega \big(\partial_t(|z|^2)  \left( \varphi'' - \Delta \varphi \right)\big)(T)  dx
+\frac{s}{2} \int_\Omega \big(|z|^2 \left( \partial_t^2 - \Delta \right) \varphi'\big)(T) dx
\\
&\quad -s \int_\Omega \big( | \nabla z |^2 \varphi'\big)(T) dx -s^3 \int_\Omega \big(|z|^2  \varphi' \left( |\varphi'|^2 - |\nabla \varphi|^2 \right)\big)(T) dx.
\end{split}
\end{equation}
The boundary term \(\mathcal{B}\) defined by (\ref{B}) becomes
\begin{equation}\label{B0}
\begin{split}
- \mathcal{B}  \leq C \int_{\Sigma_0} s\lambda\varphi |\partial_{\nu} z|^2 \partial_{\nu} (|x-x_0|^2) ds dt
\leq C \int_{\Sigma_0} s\lambda\varphi |\partial_{\nu} z|^2  ds dt.
\end{split}
\end{equation}
Furthermore, we have
\begin{equation}\label{J4}
\begin{split}
&\left|J_4+ 2(1 +\beta) \int_\Omega \big(s\lambda\varphi z'\big)(T) dx -(1 +\beta) \int_\Omega \big((s\lambda\varphi)' |z|^2\big)(T) dx \right|
\\
& \quad \leq C \int_{\Omega} \left( s|\nabla z(x,T)|^2 +s|\partial_t z(x,T)|^2+ s^3 |z(x,T)|^2 \right)
\end{split}
\end{equation}
Combining (\ref{J4}), (\ref{B0}) and (\ref{Ps}), we obtain
\begin{equation}\label{Ps0}
\begin{split}
 \int_Q s &\left( |\nabla z|^2 + |z'|^2 + s^2 |z|^2 \right) dx dt
\\
&\leq C\Big(\| P_s z \|^2 +\int_{\Sigma_0} s|\partial_{\nu} z|^2  ds dt+\int_{\Omega} ( s|\nabla_{x,t} z(x,T)|^2 + s^3 |z(x,T)|^2 )dx\Big).
\end{split}
\end{equation}
By substituting  $z=e^{s\varphi}v$  into (\ref{Ps0}), we prove Lemma \ref{lem}.

\end{proof}


\begin{thebibliography}{100}

\bibitem{BK1981}A. Bukhgeimand M. Klibanov, Global uniqueness of a class of multidimensional inverse
 problems, Sov. Math. Dokl., 24(1981), 244–247.

\bibitem{HIY} X. Huang,, O. Imanuvilov and M. Yamamoto,  Stability for inverse source problems by Carleman estimates. Inverse Problems, 36(2020), 125006.

\bibitem{Beilina2012}
L. Beilina L and M. Klibanov, 
Approximate Global Convergence and Adaptivity for Coefficient Inverse Problems,
Berlin: Springer; 2012.

\bibitem{Bellassoued2017}
M. Bellassoued M and M. Yamamoto M ,
Carleman Estimates and Applications to Inverse Problems for Hyperbolic Systems,
Tokyo: Springer-Japan; 2017.




\bibitem{Cannarsa2019a}
P. Cannarsa P, G. Floridia, F. Gölgeleyen and M. Yamamoto, 2019
Inverse coefficient problems for a transport equation by local Carleman estimate,
Inverse Problems, 35(2019), 105013.



\bibitem{Fu2019}
X. Fu, Q. Lü and X. Zhang, Carleman Estimates for Second Order Partial Differential Operators and Applications,
Berlin: Springer; 2019.

\bibitem{Gölgeleyen2016}
F. Gölgeleyen F and M. Yamamoto,
Stability for some inverse problems for transport equations, SIAM J. Math. Anal., 48(2016), 2319-2344.



\bibitem{Imanuvilov1998}
O. Imanuvilov O and M. Yamamoto, 
Lipschitz stability in inverse parabolic problems by the Carleman estimate,
Inverse Problems, 14(1998), 1229-1245.

\bibitem{Imanuvilov2001}
O. Imanuvilov O and M. Yamamoto,
Global Lipschitz stability in an inverse hyperbolic problem by interior observations, Inverse Problems, 17(2001), 717-728.

\bibitem{Imanuvilov2001a}
O. Imanuvilov O and M. Yamamoto, 
Global uniqueness and stability in determining coefficients of wave equations, Commun. Part. Differ. Equ., 26(2001), 1409-1425.

\bibitem{Imanuvilov2003}
O. Imanuvilov O and M. Yamamoto, 2003
Determination of a coefficient in an acoustic equation with a single measurement, Inverse Problems, 19(2003), 157-171.

\bibitem{Isakov1990}
V. Isakov,
Inverse Source Problems,
 Providence (RI): American Mathematical Society; 1990.

\bibitem{JLY2017} D. Jiang, Y. Liu and M. Yamamoto, Inverse source problem for the hyperbolic equation
with a time-dependent principal part, J. Differ. Equ, 262 (2017), 653–681.

\bibitem{K2002} M. V. Klibanov, Carleman estimates and inverse problems: uniqueness and convexification
of multiextremal objective functions, (2002), 219–252.


\bibitem{Y1995} M. Yamamoto, Stability, reconstruction formula and regularization for an inverse source
hyperbolic problem by a control method, Inverse Probl, 11 (1995), 481–496.

\bibitem{Y1999} M. Yamamoto, Uniqueness and stability in multidimensional hyperbolic inverse problems,
J. Math. Pures Appl., 78 (1999), 65–98.




\end{thebibliography}
\end{document}